\documentclass[11 pt]{amsart}
\usepackage{amssymb}
\usepackage{amsmath}
\usepackage{amsfonts}
\usepackage{graphicx}
\usepackage{amsthm}
\usepackage{enumerate}
\usepackage[mathscr]{eucal}
\usepackage{verbatim}
 \theoremstyle{plain}
\newtheorem{theorem}{Theorem}[section]
\newtheorem{lemma}[theorem]{Lemma}
\newtheorem{proposition}[theorem]{Proposition}
\newtheorem{corollary}[theorem]{Corollary}

\numberwithin{equation}{section}
\theoremstyle{plain}

\numberwithin{equation}{section}
\theoremstyle{remark}

\def\defeq{:=}

\def\intslash{\rlap{\kern  .32em $\mspace {.5mu}\backslash$ }\int}
\def\qsl{{\rlap{\kern  .32em $\mspace {.5mu}\backslash$ }\int_{Q_x}}}

\def\vth{\vartheta}

\def\emph#1{{\it #1 }}

\def\inn#1#2{\langle#1,#2\rangle}

\def\lc{\lesssim}

\def\ga{\gamma}

\def\ka{\kappa}

\def\la{\lambda}

\def\fA{{\mathfrak {A}}}

\def\fS{{\mathfrak {S}}}

\def\fc{{\mathfrak {c}}}

\def\fv{{\mathfrak {v}}}

\def\bbR{{\mathbb {R}}}

\def\bbZ{{\mathbb {Z}}}
\def\cA{{\mathcal {A}}}
\def\cB{{\mathcal {B}}}
\def\cC{{\mathcal {C}}}

\def\cI{{\mathcal {I}}}
\def\cJ{{\mathcal {J}}}
\def\cK{{\mathcal {K}}}

\def\cR{{\mathcal {R}}}

\def\cU{{\mathcal {U}}}

\def\C{{\hbox{\bf C}}}

 at 10 true pt

\def\be#1{\begin{equation}\label{#1}}
\def\ee{\end{equation}}
\def\bas{\begin{align*}}
\def\eas{\end{align*}}
\def\bi{\begin{itemize}}
\def\ei{\end{itemize}}

\def\emph#1{{\it #1}}
\def\textbf#1{{\bf #1}}
\def\intslash{\rlap{\kern  .32em $\mspace {.5mu}\backslash$ }\int}
\def\qsl{{\rlap{\kern  .32em $\mspace {.5mu}\backslash$ }\int_{Q_x}}}

\begin{document}


\title
[Restriction of Fourier transforms to curves, II]
{Restriction of Fourier transforms to curves\\
 II:
Some classes with vanishing torsion}

\author[]
{Jong-Guk Bak \ \ Daniel M. Oberlin \  \  Andreas Seeger}

\address {J. Bak \\ Department of Mathematics
and the Pohang Mathematics
Institute \\ Pohang University of Science and Technology
\\
Pohang 790-784, Korea}
\email{bak@postech.ac.kr}

\address
{D. M.  Oberlin \\
Department of Mathematics \\ Florida State University \\
 Tallahassee, FL 32306}
\email{oberlin@math.fsu.edu}

\address{A. Seeger   \\
Department of Mathematics\\ University of Wisconsin-Madison\\Madison,
WI 53706, USA}
\email{seeger@math.wisc.edu}


\thanks{J.B. was supported in part by
 grant
R01-2004-000-10055-0 of the Korea Science and Engineering Foundation.
D.O. was supported in part by NSF grant DMS-0552041.
A.S. was supported in part by NSF grant DMS-0200186.
}

\begin{abstract}
We consider  the Fourier restriction operators
associated to certain degenerate curves in $\mathbb R^d$
for which the highest torsion vanishes.  We prove estimates with respect to
affine arclength and with respect to the Euclidean arclength measure on the curve.
The estimates have certain uniform features, and the affine arclength results cover families of flat curves.
\end{abstract}

\maketitle

\section{Introduction}\label{intro}

We suppose that $\gamma$ is a curve in $\mathbb R ^d$ and consider the
problem of obtaining $L^p \rightarrow L^q$ bounds for the restriction of
the Fourier transform to $\gamma$. This problem has a long and
interesting
history which is described at length in \cite{D2} and \cite{BOS}. Though we will not
repeat much of that description here, we recall one of the main results
from \cite{BOS}, concerning the moment curve
$\gamma_0 (t)=(t ,t^2 ,\dots ,t^d )$ in dimension $d\ge 3$.
Write $p_d =\frac{d^2 +d+2}{d^2 +d}$. Then there is the restricted 
strong type 
inequality
\begin{equation}\label{moment}
\Bigl(\int_a^b |\widehat f (\gamma (t))|^{p_d}dt \Bigr)^{1/p_d}
\le C(\gamma )\,\|f\|_{L^{p_d ,1}(\mathbb R ^d )},
\end{equation}
for all  Schwartz functions $f$ on $\mathbb R ^d$. The  estimate
\eqref{moment} is, as described in \cite{BOS}, best possible and
yields all other $L^p \rightarrow L^q$ restriction results for the
moment curve $\gamma_0$ by interpolation with the trivial 
$L^1\to L^\infty$
estimate. It is natural to wonder what happens to \eqref{moment}
when $\gamma _0$ is replaced by more general curves. If $\gamma
:[a,b]\rightarrow \mathbb R ^d$ is nondegenerate in the sense that
for each $t\in [a,b]$ the derivatives $\gamma '(t),\gamma ''
(t),\dots ,\gamma ^{(d)}(t)$ are linearly independent, then the
analogue of \eqref{moment} is proved in \cite{BOS}. But if one
attempts to go further by dropping the hypothesis of
nondegeneracy, it is easy to see that the exact analogues of
\eqref{moment} and its interpolants may fail. There are then two
possibilities which have been considered in the literature. The
first is to \lq\lq dampen" the measure $dt$ by introducing a
weight $w(t)$ which is small where $\gamma$ is degenerate, to
replace $dt$ with $w(t)\,dt$, and then to attempt to obtain the
exact analogue of \eqref{moment}. The second approach is to retain
$dt$ for the reference measure and to see what changes must then
be made in order to obtain sharp restriction results. In this
paper we explore both approaches, but only for $\gamma$ of the
form
\begin{equation}\label{specialga}
\gamma (t)=\Bigl( t,\frac{t^2}{2},\dots ,\frac{t^{d-1}}{(d-1)!},\phi (t)
\Bigr) .
\end{equation}
These curves are termed
{\it simple}
in \cite{DM1} and are distinguished
by the fact that only the highest torsion may vanish.

Concerning the first approach, it was observed in \cite{DM1} that if
$\gamma$ is as in \eqref{specialga}, then the
correct weight $w(t)$  is given by
\begin{equation}\label{correctweight}
w(t)=|\phi ^{(d)}(t)|^{\frac{2}{d(d+1)}}.
\end{equation}
Then the measure $w(t)\,dt$ is, up to a constant depending only on the
dimension, the affine arclength measure on $\gamma$. Here we
have the following result.

\begin{theorem} \label{flatthm}
Fix $d\ge 2$. Suppose $0 \le a<b\le\infty$ and let $\gamma$ be of
the form \eqref{specialga} where $\phi$ is a $C^d$ function on
$(a,b)$ for which the derivatives $\phi ', \dots ,\phi ^{(d)}$ are
nonnegative and nondecreasing on $(a,b)$ and for which $\phi
^{(d)}$ satisfies the condition
\begin{equation}\label{techcondnew}
\Big(\prod_{j=1}^d \phi^{(d)}(s_j)\Big)^{1/d} \le A \,  \phi^{(d)}
\big(\tfrac{s_1+\dots+s_d}{d}\big)
\end{equation}
for all $s=(s_1,\dots, s_d)$ with $a<s_1\le s_2\le \dots\le s_d < b$.

Suppose $1\le P <\tfrac{d^2+d+2}{d^2+d}$, and $1-\frac{1}{P}=
\frac{2}{d(d+1)}\frac{1}{Q}$.
Then there is $C(d,P)<\infty$ so that for all $g\in
L^{P}(\mathbb R ^d)$
\begin{equation}\label{affrestrflat}
\Big(\int_a^b
 |\widehat g (\gamma(t))|^{Q} w(t)\,  dt\Big)^{1/Q} \le
C(d,P)\, A^{1-1/P}\,\|g\|_{L^{P}(\mathbb R ^d)}.
\end{equation}
\end{theorem}
\medskip

\noindent The proof of Theorem \ref{flatthm} is analogous to the
proof of Theorem 1.3 in \cite{BOS}. The range of indices in
Theorem \ref{flatthm} is the range given by interpolating  the
$L^{p_d,1}\to L^{p_d}$ estimate \eqref{moment} with the trivial
$L^1\to L^\infty$ estimate, and it would be interesting to know if
the endpoint result (the exact analogue of  \eqref{moment}) holds
for the curves of Theorem \ref{flatthm}.

In the case $d=2$ it follows from
\cite{sj} that the conclusion of Theorem \ref{flatthm} holds with $A=1$
(and without any additional hypotheses like \eqref{techcondnew}).
For many interesting examples a slightly stronger condition holds where the
arithmetic mean in the argument of $\phi^{(d)}$ on the right hand side of \eqref{techcondnew}
is replaced by a geometric mean, {\it i.e.},
\begin{equation}\label{techcond}
\Big(\prod_{j=1}^d \phi^{(d)}(s_j)\Big)^{1/d} \le A \,  \phi^{(d)}
\big(\root \uproot 1d \of{s_1 \cdots s_d}\big).
\end{equation}
It is obvious that condition \eqref{techcond}
holds  for $\phi(t)=t^\beta$, $\beta\ge d$ on the interval
$(0,\infty)$; in particular \eqref{techcond}
is satisfied with $A=1$. Moreover, if for $t\ge 0$ we define
$\phi_0(t)=t^{\beta}$ for some $\beta >d$,  and for $n\ge 1$,
$$
\phi_n(t)=
\int_0^t {(t-u)^{d-1}}
\exp\Big(-\tfrac{1}{\phi_{n-1}^{(d)}(u)}\Big) du,
$$
then
$\phi_n$ satisfy \eqref{techcond} with $A=1$ on
$(0,\infty)$ (see \S \ref{flex}). This
yields  a  sequence of functions which are
progressively flatter at the origin
for which the restriction theorem holds uniformly (i.e., with constant
depending only on the Lebesgue space indices).
These two observations raise the interesting question of whether or not
the hypothesis \eqref{techcond} in Theorem \ref{flatthm} can
be dropped to yield, subject to $\phi$'s being sufficiently monotone,
a uniform restriction theorem for the curves \eqref{specialga}.

Regarding the second of the above-mentioned possibilities, keeping
the measure $dt$, Drury and Marshall \cite{DM2} proved sharp
results for classes of finite type curves. Here we are aiming for
a result for curves of the form \eqref{specialga} which is
expressed in terms of a natural geometric condition and also has a
certain uniform feature.

We will say that a set $E$ in $\mathbb R ^d$ is a {\it
parallelepiped } if $E$ is a translate of a set of the form
$\{\sum_{j=1}^d t_j x_j :0\le t_j \le 1\}$ where the  $x_j
\in\mathbb R ^d$ are linearly independent. Given $\gamma$ we shall
write $\lambda_\ga$ for the measure on $\gamma$ given by
$$\inn{d\lambda_\ga}{f} =\int f(\gamma(t))dt.$$ We denote Lebesgue
measure in $\bbR^d$ by $m_d$.

\begin{theorem} \label{euclthm}
Suppose $-\infty <a<b<\infty$ and let
$\gamma$ be of the form \eqref{specialga} where $\phi$ is a $C^d$ function
on $(a,b)$ for which the derivatives $\phi ' ,\dots ,\phi^{(d)}$ are
nonnegative and
nondecreasing on $(a,b)$. Suppose that $\alpha \in (0,\tfrac{2}{d(d+1)}]$
if $d\geq
3$ and that $\alpha \in (0,1/3)$ if $d=2$.
Suppose also  that the estimate
\begin{equation}\label{hypothesis2}
\lambda_\ga (E)\leq B\, m_d (E)^{\alpha}
\end{equation}
holds for some $B>0$ and for all parallelepipeds $E\subset\mathbb R^d$.

Then there is $C(d,\alpha )<\infty$ so that for all
$g\in L^{1+\alpha,1}(\mathbb R ^d )$
\begin{equation}\label{conclusion2}
\Big(\int_a^b
|\widehat g (\gamma(t))|^{1+\alpha}  dt\Big)^{1/(1+\alpha )} \le
C(d,\alpha )\, B^{\frac{1}{1+\alpha} }\,\|g\|_{L^{1+\alpha ,1}(\mathbb R
^d)}.
\end{equation}
On the other hand, if the estimate
\begin{equation}\label{hypothesis3}
\Big(\int_a^b
 |\widehat g (\gamma(t))|^{Q}  dt\Big)^{1/Q} \le
 c^{{1}/{Q} }\,\|g\|_{L^{P}(\mathbb R^d)}
\end{equation}
holds for some $P$ and $Q$ satisfying
$1-\frac{1}{P}=\frac{\alpha}{Q}$, then \eqref{hypothesis2}
holds for all parallelepipeds $E$ with $B$ replaced by $C(d,p)\, c$.
\end{theorem}
\noindent The proof of Theorem \ref{euclthm} is analogous to the proof
of  \eqref{moment} given in \cite{BOS}. Interpolation of
\eqref{conclusion2} with the trivial $L^1$ estimate yields the
estimate
\eqref{hypothesis3} whenever $1\le P <1+\alpha$ and
$1/P'=\alpha /Q$. It would be interesting to know whether in the generality of Theorem \ref{euclthm} the exponent
$1+\alpha$ is sharp when $\alpha <2/(d^2 +d)$ or whether there is always
$P(\alpha )>1+\alpha$
such that \eqref{hypothesis2} implies \eqref{hypothesis3} whenever $1\le
p <P(\alpha )$ and $1/P'=\alpha /Q$. For many concrete examples such
improvements can  indeed be obtained by
rescaling arguments from the nondegenerate case --
for this and related observations see \S\ref{further}.

{\it This paper:}
In order to prove Theorem \ref{flatthm}  we shall use the method of offspring curves that originated in \cite{D}, and was further developed in \cite{DM1},
\cite{DM2} and \cite{BOS}.  The  crucial
 technical point is to give lower bounds for a certain
Jacobian of a change of variable, estimate \eqref{Jlow} below.
 The new features about Theorem \ref{flatthm}
concern the verification of this estimate, and the technical
details are contained in \S \ref{flat}.
The proof of Theorem \ref{flatthm} is then discussed in \S\ref{pfofthm}
(a reader not familiar with the method should start reading here).
In \S\ref{examples} we  discuss some examples to which Theorem
\ref{flatthm} can be applied. Sections
\S\ref{eucl} and \S\ref{lemmas}  contain the proof of  Theorem \ref{euclthm}.
In \S\ref{further} we show how
Theorem \ref{euclthm} can be extended for some classes of examples.

\section{The main technical estimate}\label{flat}
In this section we assume that $\phi$ is defined on $[a,b]$,
$0\le a<b$ and assume that the derivatives of $\phi$ up to order $d$
are positive and nondecreasing on $(a,b)$.

We establish some notation. For a vector $x\in \mathbb R^d$
let
$V_d(x)$ be the determinant of the
$d\times d$ Vandermonde matrix:
\begin{equation}
\label{vanderm}
 V_d(x)= \prod_{1\le i<j\le d} (x_j-x_i).
\end{equation}
For $h=(h_1,\dots, h_{d-1})\in (\Bbb R_+)^{d-1}$ define
$\kappa(h)\in [0,\infty)^d$ by
\begin{equation*}
\kappa_1(h)=0, \quad \kappa_{j}(h)= h_1+\dots+ h_{j-1}, \quad 2\le j\le d
\end{equation*}
and put
\begin{equation*}
v(h)\equiv v_d(h)=V_d(\ka(h)).
\end{equation*}
If $\gamma:(a,b)\rightarrow \mathbb R ^d$ and if $a<t<b-\ka_d (h)$,
we write
\begin{equation}
\label{Gath}
\Gamma(t,h)=\sum_{j=1}^d \gamma(t+\ka_j(h)).
\end{equation}
Following the
 terminology of Drury and Marshall \cite{DM1} we  call   $\Gamma (\cdot ,h)$,
for fixed $h$,  an {\it offspring curve} of
$\gamma$.



Denote by
 $J_\phi(t,h) $ the
Jacobi-deter\-minant of the transformation
$(t,h)\mapsto \Gamma(t,h);$
that is
\begin{equation}\label{Jphidet}
J_\phi(t,h)=\det
\begin{pmatrix} \frac{\partial\Gamma}{\partial t}
& \frac{\partial\Gamma}{\partial h_1}
&\dots
& \frac{\partial\Gamma}{\partial h_{d-1}}
\end{pmatrix}.
\end{equation}

As in \cite{DM1}  it will be crucial to verify the identity
\begin{equation}
\label{Jlow}
|J_\phi(t,h)|\ge \sigma v(h)
\Big(\prod_{i=1}^d \phi^{(d)}(t+\kappa_i(h))\Big)^{1/d}.
\end{equation}
Here we prove
\begin{proposition}\label{jacobicond} Let $0\le a<b\le 1$. Suppose that
$\phi ^{(d)}$ is nonnegative
and nondecreasing on $(a,b)$, and  that
for   any $a<s_1 \leq\cdots
\leq s_d\le b$,
the condition
\eqref{techcondnew} is satisfied. Then
condition \eqref{Jlow} holds with $$\sigma=c_0(d)A^{-1}$$  for all
$(t,h)$ such that $a\le t\le b$, $h\in (0,b)^{d-1}$, and
$t+\ka_d(h)\le b$.
\end{proposition}

The proof of Proposition \ref{jacobicond} uses the following
technical lemma.

\begin{lemma}\label{lin} Fix $\lambda \in (0,1)$. Suppose
$$
a_1 <b_1 \leq a_2 <b_2 \leq \cdots \leq a_N <b_N .
$$
Suppose that, for $m=1,\dots ,M$,  $l_m$ is a function of $t=(t_1 ,
\dots ,t_N )$ having one of the three following forms:
$$
l_m (t)=
\begin{cases} t_k -t_j \ \ \text{  for  some}\ \ 1\leq j<k\leq N, \
\text{ or}
\\
d_j -t_j \ \ \text{ for  some}\ \ d_j \geq b_j , \ \text{ or}
\\
t_j -c_j \ \ \text{ for  some}\ \ c_j \leq a_j .
\end{cases}
$$
Suppose that
$\lambda _j \in (0,1)$ and $\lambda _j
\leq \lambda$, for $j=1,\dots ,N$.
Let $\cR_N(a,b,\la)$ be the region of all
$t=(t_1,\dots, t_N)\in \bbR^N$ satisfying
$  (1-\la _j )a_j +\la _j b_j
\le t_j \le b_j$ for $j=1,\dots, N$.
Then
\begin{multline*}
\int\limits_{
\cR_N(a,b,\la)}
\prod\limits_{m=1}^M l_m (t) \, dt_N \cdots dt_1
\\
\geq C(M,\lambda )^N
\int\nolimits _{a_1}^{b_1}\cdots \int\nolimits_{a_N}^{b_N}
\prod\limits_{m=1}^M l_m (t) \,  dt_N \cdots dt_1 .
\end{multline*}
\end{lemma}
\begin{proof}[Proof of Lemma \ref{lin}] An easy induction argument shows
that
it is enough to prove the lemma when $N=1$. A translation and then
a scaling
reduce that case to the inequality
\begin{equation}\label{Nequal1}
\int\nolimits_{1}^{1/\la}\prod\limits_{m=1}^M l_m (t) \, dt
\geq C(M,\lambda )
\int\nolimits_{0}^{1/\la}\prod\limits_{m=1}^M l_m (t) \, dt
\end{equation}
where
$$
l_m (t)=\begin{cases}
d_m -t \,&\text{ for  some }\ d_m \geq 1/\la,  \ \text{ or}\\
t -c_m \,&\text{ for  some }\ c_m \leq 0.
\end{cases}
$$
It is clear that
\eqref{Nequal1}
is true when $M=0$. So assume that \eqref{Nequal1}  is true for $M-1$.
Suppose first that at least one of the functions $l_m$ is increasing,
say $l_M (t)=t-c$. Then, by the inductive assumption,
$$
\int\nolimits_{1}^{1/\la}
\prod\limits_{m=1}^{M-1} l_m (t) \, dt
\geq \frac{C(M-1,\lambda )}{1-C(M-1,\lambda )}
\int\nolimits_{0}^{1}\prod\limits_{m=1}^{M-1} l_m
(t) \, dt.
$$
Thus
\begin{multline*}
\int\nolimits_{1}^{1/\la}\prod\limits_{m=1}^{M-1} l_m
(t)(t-c) \, dt \geq (1-c)\int\nolimits_{1}^{1/\la}\prod\limits_{m=1}^{M-1}
l_m (t) \, dt
\\ \geq
\frac{C(M-1,\lambda )}{1-C(M-1,\lambda )}
\int\nolimits_0^{1}\prod\limits_{m=1}^{M-1} l_m
(t)(t-c) \, dt,
\end{multline*}
and this is equivalent to \eqref{Nequal1}  with $C(M-1 ,\lambda )$ instead
of $C(M ,\lambda )$. There\-fore we can assume that $l_m (t)=d_m -t$
for all $m$. There are now two cases to consider. First suppose that one
of the $d_m$'s, say $d_M$, exceeds $2/\la$.
Let $\tau=(1+1/\la )/2$. Then
\begin{multline}\label{so-est2}
\int\nolimits_{0}^{1}\prod\limits_{m=1}^{M} l_m (t) \, dt
\leq d_M
\int\nolimits_{0}^{1}\prod\limits_{m=1}^{M-1} l_m (t) \, dt
\\
\leq
d_M \frac{1-C(M-1,\lambda )}{C(M-1,\lambda )}
\int\nolimits_{1}^{1/\la} \prod\limits_{m=1}^{M-1} l_m (t) \, dt.
\end{multline}
We further estimate
\begin{equation}\begin{aligned}
\label{so-est3}
\int\nolimits_{1}^{1/\la} \prod\limits_{m=1}^{M-1} l_m (t) \, dt
&
\leq 2\int\nolimits_{1}^{\tau} \prod\limits_{m=1}^{M-1} l_m
(t)
\, dt
\\
&
\le
\frac{2}{d_M -\tau}
\int\nolimits_{1}^{\tau} \prod\limits_{m=1}^{M} l_m (t) \, dt
\\
&
\leq
\frac2{d_M -\tau}
\int\nolimits_{1}^{1/\la} \prod\limits_{m=1}^{M} l_m (t) \, dt
\end{aligned}
\end{equation}
where the first inequality follows
because $\prod_1^{M-1}l_m(t) dt$ is decreasing.
Since $d_M \le 2/\la$, we have $d_M /(d_M -\tau )\le 2$.
Combined with \eqref{so-est2} and \eqref{so-est3},
this implies \eqref{Nequal1} if one of the $d_m$'s exceeds $2/\la$.
If, on the other hand, $d_m\le 2/\la$ for all $m$, then
the crude estimates
$$
\int_0^{1/\la}\prod_{m=1}^M l_m (t)\, dt\le \frac{2^M}{\la ^{M+1}}
$$
and
$$
\int_1^{1/\la}\prod_{m=1}^M l_m (t)\, dt
\ge
\int_1^{1/\la}\Bigl(\frac{1}{\la}-t\Bigr)^M \, dt=\frac{(1/\la
-1)^{M+1}}{M+1},
$$
give \eqref{Nequal1} again and conclude the proof of Lemma \ref{lin}.
\end{proof}

It will be useful to write the Jacobian $J_\phi(\cdot,h)$ as a convolution
with a nonnegative function,
 depending on the parameter $h\in (\bbR_+)^{d-1}$.

To this end we define  for $h_1\ge 0$
\begin{equation}\label{psidef2}
\Psi_2(t;h_1)= \chi_{[0,h_1]}(t).
\end{equation}
For $d\ge 3$ and $t\le h_1+\dots+h_{d-1}$ we set
\begin{align*}
\cR_{d-1}(t,h)=\big\{ \sigma\in &\bbR^{d-1}: 0\le \sigma_1\le \min
\{h_1,t\},\\ 
&h_1+...+h_{j-1} \le \sigma_j \le h_1+...+h_j,\quad
j=2,...,d-2,
\\
&\max\{h_1+...+h_{d-2},t\} \le \sigma_{d-1}\le h_1+...+h_{d-1} \big\}\, ,
\end{align*}
and  define recursively
\begin{equation}\label{psidefd}
\Psi_d(t;h_1,\dots, h_{d-1})=
\int_{\cR_{d-1}(t,h)} \Psi_{d-1}(t-\sigma_1;\sigma_2,\dots, \sigma_{d-1})
d\sigma_1\dots d\sigma_{d-1}
\end{equation}
if $t\le h_1+\dots+h_{d-1}$; we also set
$\Psi_d(t;h)=0$   if $t\ge h_1+\dots+h_{d-1}$.

\begin{lemma}\label{Jdsphilemma}
 Let $\Psi_d$ be as in \eqref{psidef2}, \eqref{psidefd} and
let, for $s\in \bbR^d$ with
$s_1\le s_2\le ...\le s_d$,   $\cJ_d (s_1 ,\dots ,s_d;\phi)$ denote
the determinant of the $d\times d$ matrix with columns $(1,s_j ,\dots,
\tfrac{s_j^{d-2}}{(d-2)!},\phi '(s_j ))^T$.

Then \begin{equation}\label{integraldarstellung} \cJ_d(s;\phi)=
\int_{s_1}^{s_d} \Psi_d( u-s_1; s_2-s_1,\dots, s_d-s_{d-1})
\,\phi^{(d)}(u) du .
\end{equation}
\end{lemma}
\begin{proof}
If $d=2$ then the asserted formula holds since
$$\cJ_2(s_1,s_2;\phi)=\phi'(s_2)-\phi'(s_1)=\int_{s_1}^{s_2}\phi''(u)du
$$ and $\Psi_2(u-s_1;s_2-s_1)= \chi_{[s_1,s_2]}(u)$.

We now argue by induction and assume $d\ge 3$.

We first note
by expanding
$\partial_1\dots\partial_{d-1} \cJ_d$ with respect to the last column
that
$$\partial_{s_1}\dots\partial_{s_{d-1}} \cJ_d(s_1,\dots,s_d;\phi)=
(-1)^{d+1}\cJ_{d-1}(s_1,\dots,s_{d-1}; \phi').
$$
Next observe that $\cJ_d(s;\phi)=0$  if $s_1=s_2$ and that
$\partial_1\dots\partial_k \cJ_d(s;\phi)=0$ if $s_{k+1}=s_{k+2}$ and
$k\le d-2$. Thus we repeatedly integrate and see that

\begin{align}\label{firstintegration}
&\quad\cJ_d(s;\phi)
\\
&=(-1)^{d-1}\int_{s_1}^{s_2}\dots\int_{s_{d-1}}^{s_d }
\partial_{s_1}\dots\partial_{s_{d-1}} \cJ_d(\sigma_1,\dots,\sigma_{d-1},s_d;\phi)
\,d\sigma_{d-1}\dots d\sigma_1
\notag
\\&=
\int_{s_1}^{s_2}\dots\int_{s_{d-1}}^{s_d }
\cJ_{d-1}(\sigma_{1},\dots,\sigma_{d-1};\phi')\,
d\sigma_{d-1}\dots d\sigma_{1}.
\notag\end{align}

Thus by the induction  hypothesis
\begin{multline} \label{cJint}
\cJ_d(s;\phi)=
\int_{s_1}^{s_2}\dots\int_{s_{d-1}}^{s_d }\int_{\sigma_1}^{\sigma_{d-1}}
\phi^{(d)}(u) \quad\times\\
\Psi_{d-1}(u-\sigma_1 ;
\sigma_2-\sigma_1,\dots, \sigma_{d-1}-\sigma_{d-2})
\, du \,  d\sigma_{d-1}\dots d\sigma_{1}
\end{multline}
and  by Fubini's theorem
this can be written in the form
\begin{multline*}
\cJ_d(s;\phi)=\int_{s_1}^{s_d} \phi^{(d)}(u)
\int_{\tau\in \Omega(u)}
\Psi_{d-1}(u-\tau_1 ;
\tau_2-\tau_1,\dots, \tau_{d-1}-\tau_{d-2}) d\tau du
\end{multline*}
where $\Omega(u)$ consists of those $\tau \in \bbR^{d-1}$ for which
$s_i\le \tau_i\le s_{i+1}$, $i=1,\dots d-1$ and $\tau_1\le u\le \tau_{d-1}$.

We change variables $\tau_i=s_1 +\sigma_i$ for $i=1,...,d-1,$
so that
$\tau\in \Omega(u)$ corresponds to
$\sigma\in \cR_{d-1}(u-s_1,h)$ with $h_i=s_{i+1}-s_i$. Thus from
the definition \eqref{psidefd} we obtain
\begin{equation*}
\cJ_d(s;\phi)=
\int_{s_1}^{s_d} \Psi_d(
u-s_1 ; s_2-s_1,\dots, s_d-s_{d-1})
\,\phi^{(d)}(u) du
\end{equation*}
which yields the assertion.
\end{proof}

\begin{lemma}\label{Jintegral}
Let $\Psi_d$ as in \eqref{psidef2}, \eqref{psidefd} and let
$$g_d(t,h)=t+\frac 1d\sum_{i=1}^d \kappa_i(h).$$
Then  $\Psi_d$ satisfies
\begin{equation}\label{intpsilow}
\int^{t+\ka_d(h)}_{g_d(t,h)}
\Psi_d(u-t;h) du\ge c(d) v(h)
\end{equation}
where
 $c(d)>0$.
\end{lemma}
\begin{proof}
First, in order to prepare for  the proof of \eqref{intpsilow}, we
 observe that \eqref{firstintegration} for  the special case
$\phi(s)=s^d/d!$ gives us the formula for the Vandermonde determinant
$V_d(s)= \prod_{j=1}^{d-1} (j!) \cJ_d(s,\phi)$ in all dimensions namely
\begin{equation} \label{vandintegration}
\begin{aligned}
V_n(s_1,\dots, s_d) = (n-1)!\int_{s_1}^{s_2}\dots\int_{s_{n-1}}^{s_n }
V_{n-1}(\sigma_1,\dots,\sigma_{n-1})
d\sigma_{n-1}\dots d\sigma_{1}.
\end{aligned}
\end{equation}
We  now use  Lemma \ref{lin}
to establish the following inequality
for all $n\ge 2$.
Suppose that $0\leq a_1 \leq \cdots \leq a_n$ and let
$$\cU_{n-1}(a)=\{x'=(x_1,\dots,x_{n-1})\in (\Bbb R_+)^{n-1}:
\frac{1}{n-1}\sum_{i=1}^{n-1} x_i\ge
\frac{1}{n}\sum_{k=1}^{n} a_k
\}.$$

Then
\begin{multline} \label{Vmixedbd}
\int\nolimits_{a_1}^{a_2}\cdots\int\nolimits_{a_{n-1}}^{a_n}
V_{n-1}(x')\cdot\chi_{\cU_{n-1}(a)}(x')\, dx' \geq
C(n)\,V_n(a_1 ,\dots ,a_{n}).
\end{multline}

To check \eqref{Vmixedbd}, note that if $\la_j =(n-j)/n,$ then
the left hand side of \eqref{Vmixedbd} certainly exceeds
$$
\int\nolimits_{\la_1 a_1 +(1-\la_1 )a_2}^{a_2}
\cdots \int\nolimits_{\la_{n-1} a_{n-1} +(1-\la_{n-1} )a_n}^{a_n}
V_{n-1}(x)\, dx_{n-1}\cdots dx_{1} .
$$
By Lemma \ref{lin} this expression is bounded below by a positive
constant times the integral of $V_{n-1}$ over the entire rectangle
$\prod_{i=1}^{n-1} [a_i,a_{i+1}]$, and by
\eqref{vandintegration},
that integral is equal to $C(n)\, V_{n}(a_1,\dots, a_n)$.

We shall now prove \eqref{intpsilow}.
The case $d=2$ is immediate since
$\Psi_2 (\cdot \,;h)=\chi _{[0,h_1 ]}$ and $v(h_1 )=h_1$: we find that
\eqref{intpsilow} holds with $c(2)=1/2$.
Now we argue by induction and assume that \eqref{intpsilow} holds if
$d-1\ge2$. With $s_j =t+\kappa _j (h)$ we use
\eqref{cJint}
for a $\phi$ with  $\phi^{(d)}(u)=1$ for $u\ge \overline s
=(s_1 +\cdots +s_d )/d$
and
$\phi^{(d)}(u)=0$ for $u< \overline s$.
We thereby obtain
\begin{equation*}
\begin{aligned}
&\int_{g_d(t,h)}^{t+\kappa_d (h)} \Psi_d(u-t;h) \, du
=\cJ_d (s;\phi )
\\
&\,=\int_{s_1}^{s_2}\dots\int_{s_{d-1}}^{s_d
}\int_{\sigma_1}^{\sigma_{d-1}}
\,\,\,\,
\chi_{\{u\ge \overline s \}}(u)\,\,\,\,\times
\\
&\qquad\qquad
\Psi_{d-1}(u-\sigma_1 ;
\sigma_2-\sigma_1,\dots, \sigma_{d-1}-\sigma_{d-2})\,
du \,  d\sigma_{d-1}\dots d\sigma_{1}
\\
&\,\ge\,\int\nolimits_{\la_1 s_1 +(1-\la_1 )s_2}^{s_2}
\cdots
\int\nolimits_{\la_{d-1} s_{d-1} +(1-\la_{d-1} )s_d}^{s_d}
\int_{\sigma_1}^{\sigma_{d-1}}
 \,\,
\chi_{\{u\ge \overline \sigma \}}(u)\,\,\times
\\
&\qquad\qquad\Psi_{d-1}(u-\sigma_1 ;
\sigma_2-\sigma_1,\dots, \sigma_{d-1}-\sigma_{d-2})\,
du \,  d\sigma_{d-1}\cdots d\sigma_{1},
\end{aligned}
\end{equation*}
where $\la _j =(d-j)/d$.
Here the inequality follows because the conditions $\sigma_j\ge\la_j s_j
+(1-\la _j )s_{j+1}$
and $u\ge \overline \sigma =(\sigma_1 +\cdots +\sigma _{d-1})/(d-1)$
together imply
$u\ge \overline s$.
It follows from the induction hypothesis that
\begin{multline*}
\int_{\sigma_1}^{\sigma_{d-1}}
\chi_{\{u\ge  \overline \sigma \}}(u)
\Psi_{d-1}(u-\sigma_1 ;
\sigma_2-\sigma_1,\dots, \sigma_{d-1}-\sigma_{d-2})\,
du \, \\
\ge c(d-1) V_{d-1}(\sigma_1,\sigma_2,\dots, \sigma_{d-1}).
\end{multline*}
Therefore
\begin{multline*}
\int_{g_d(t,h)}^{t+\kappa _d (h)} \Psi_d(u-t;h)  du
\ge
c(d-1)\times\\
\int\nolimits_{\la_1 s_1 +(1-\la_1 )s_2}^{s_2}
\cdots
\int\nolimits_{\la_{d-1} s_{d-1} +(1-\la_{d-1} )s_d}^{s_d}
V_{d-1}(\sigma _1 ,\dots ,\sigma _{d-1})\, d\sigma _1
\cdots d\sigma _{d-1}.
\end{multline*}
With Lemma \ref{lin} and \eqref{vandintegration}, this yields
\eqref{intpsilow}.
\end{proof}

\begin{proof}[Proof of Proposition \ref{jacobicond}, conclusion]
We first observe that
\begin{equation}\label{twoJs}
J_\phi(t,h)= \cJ_d(t,t+\ka_2(h),\dots, t+\ka_d(h);\phi).
\end{equation}

Recall $g_d(t,h):=
\sum_{i=1}^d(t+\ka_i(h))/d$
so that
$t\le g_d(t,h)\le t+\ka_d(h)$.
We apply \eqref{integraldarstellung}, \eqref{intpsilow}
to get
\begin{equation*}
\begin{aligned}
J_\phi(t, h)&\ge \int\nolimits_{\overline{t}(t,h)}^{t+\ka_d(h)}
\Psi_d (u-t;h)\ \phi ^{(d)}(u)\ du
\\&\ge  \phi ^{(d)}\big(g_d(t,h)\big)
\int\nolimits_{\overline{t}(t,h)}^{t+\ka_d(h)}
\Psi_d (u-t;h) du
\\&\ge c_d \phi ^{(d)}(g_d(t,h)) \,v(h)
\ge  c_d A^{-1}
\Big(\prod_{j=1}^d \phi^{(d)}(t+\ka_j(h))\Big)^{1/d} v(h)
\end{aligned}
\end{equation*}
where we have used that $\phi^{(d)}$ is nonnegative and nondecreasing, and
in the last estimate we have employed the
hypothesized condition \eqref{techcondnew}.
\end{proof}

\medskip

\section{Proof of Theorem \ref{flatthm}}\label{pfofthm}

We first  note that $\phi$  satisfies condition \eqref{techcond} on $(0,b)$
if, and only if the function $s\mapsto \phi(bs)$ satisfies condition
\eqref{techcond} on the interval $(0,1)$. The desired estimate is invariant under the change of variable
$$x\mapsto (b^{-1}x_1,b^{-2} x_2,\dots, b^{1-d}x_{d-1} , x_d)$$
 and thus we may replace  $\phi$ by $\phi(b\cdot)$. Thus we may and shall assume  \begin{equation} b\le 1
\label{ble1}
\end{equation}
in what follows.
We shall assume also  that $\phi^{(d)}(t)$ is positive and nondecreasing in $[a,b]$ and it then suffices to prove the estimate
\eqref{affrestrflat} with the interval $(0,b)$ replaced with $(a,b)$ and $b\le 1$.

Given Proposition \ref{jacobicond} the argument is  very similar to the argument in the proof of the result for monomial curves in
\cite{BOS},
based substantially on previous ideas in papers by Christ \cite{Ch},
Drury \cite{D} and Drury and Marshall \cite{DM1},
and
the exposition will be somewhat sketchy.
We aim for an estimation of an adjoint operator and thus will set
$p= Q'=Q/(Q-1)$ and $q=P'=P/(P-1)$. Thus we fix $p<q_d=\tfrac{d^2+d+2}2$ and $q=\tfrac{d^2+d}2 p'>q_d$.
We shall now  {\it assume} that the condition \eqref{Jlow} is satisfied with a positive constant $c_0$, for all $(t,h)\in [0,1]^d$ such that
$t+\kappa_d(h)\le 1$. Note that  by Proposition \ref{jacobicond}
this assumption is implied by \eqref{techcond}.

{\bf Definition.}
\textit{  Let $0\le a<b\le 1$, $0\le M<\infty$, $\sigma>0$,
and let $\cK_{a,b,M}(\sigma)$ be the
class  of all real valued functions $\phi$ defined on $[a,b]$ for which}

\textit{(i) $\phi\in C^d([a,b])$,
$\phi^{(d)}(t)\le M$ for all $t\in [a,b]$,
 $\phi, \phi',\dots, \phi^{(d)}$ are nonnegative
on $[a,b]$, and }

\textit{(ii) for all $h\in [0,1]^{d-1}$
with $\kappa_d(h)\le b-a$ the inequality
$$J_\phi(t,h)\ge \sigma v(h) 
\Big(\prod_{i=1}^d \phi^{(d)}(t+\ka_i(h))\Big)^{1/d}$$
holds for all $t$ such that $a\le t\le b-\ka_d(h)$.  }

Let $R\ge 1$, $B_R=\{x\in \bbR^d:|x|\le R\}$,
and define
\begin{multline} \label{affrestr3}
\cA\equiv \cA(R,M,\fc) :=
\sup_{\sigma\le \fc}
\frac{\sigma}{\fc}
 \quad\times
\\
\sup_{\substack{
\phi\in \cK_{a,b,M}(\sigma)\\0\le a<b \le 1
} }
\sup_{\substack{ \|g\|_{L^{q'}(\bbR^d)}\le 1
\\ {\scriptstyle  \text{\rm supp}} (g)\subset B_R } }
\Big(\int_a^b |\widehat g (t,\dots,
\tfrac{t^{d-1}}{(d-1)!},\phi(t)) |^{p'} |\phi^{(d)}(t)|^{\frac{2}{d^2+d}}
dt\Big)^{1/p'} .
\end{multline}

Clearly $ \cA(R,M, \fc)$
is finite for each $R$ and $M$, indeed in view of $b\le 1$ we have
$ \cA(R,M, \fc)\le C_d M^{1/p'} R^{d/q'}$.
 The theorem is proved if we can show that
$\cA$ only depends on $\fc, p,  d $; in fact we will prove
that
\begin{equation}\label{cAclaim}
\cA(R,M,\fc)\le C(p,d) \, \fc^{-1/q} .
\end{equation}

The restriction inequality
$$
\Big(\int_a^b |\widehat g (\gamma (t))
|^{p'} w(t)\,  dt\Big)^{1/p'}\le \frac{\fc}{\sigma}\,\cA \,
\|g\|_{L^{q'}(B_R )}
$$
with $w=|\phi^{(d)}|^{2/[d(d+1)]}$ is equivalent to the inequality
\begin{equation}\label{extension}
\|Tf \|_{L^q(B_R)}\le \frac{\fc}{\sigma}\,\cA \, \|f\|_{L^p(wdt)},
\end{equation}
where
$$T f(x) = \int_a^b f(t) w(t) e^{-i\inn{x}{\gamma(t)}} dt.$$
For fixed $h\in (\mathbb R _+ )^{d-1}$ let
\begin{equation}
\label{Hth}
H(t,h)=\prod_{i=1}^d w(t+\ka_i(h)).
\end{equation}
With $I_h=(a,b-\ka_d(h))$
and with the convention that
$\int\cdots dt$ will mean $\int_{I_h}\cdots dt$,
we write
$$S_h[F](x) = \int e^{-i \inn{x}{\Gamma(t,h)}}
F(t,h) H(t,h)
\,dt. $$
We form $d$-fold products and, with the additional convention that
$h$ integrals are extended over the region where $\kappa_d(h)\le b$,
write
$$\prod_{i=1}^d Tf_i= \sum_{\pi\in \fS_d}\int S_{h}[F^\pi] dh$$
where
$$
F^\pi(h,t) =\prod_{i=1}^d f_{\pi(i)}(t+\kappa_i(h)).
$$
The  strategy in establishing \eqref{extension}
will be to estimate the $L^{q/d}(B_R )$ norm of
$\prod_{i=1}^d Tf_i$
by estimating
the $L^{q/d}(B_R )$ norms of
$\int S_h [F^{\pi}]\, dh$.


\begin{lemma}\label{uniform}
For  every $h$ with $\kappa_d(h)\le b-a$ the inequality
\begin{multline}
\label{uniform2} \big\|  S_{h} [F H^{-\frac{d-1} d}]\big\|_{L^q(B_R)}
\\
\le d^{d/q'} \frac{\fc}{\sigma}
\cA(R d^3, M, \sigma/d)
 \big( \int |F(t,h)|^p H(t,h)^{1/d}
dt \big)^{1/p}
\end{multline}
holds for $\phi\in \cK_{a,b, M}(\sigma)$.
\end{lemma}
\begin{proof}
Set $\overline h=d^{-1} \sum_{k=1}^{d-1}(d-k)h_k$.
A quick computation involving expansions of powers of $t$ about the point
$t+\overline h$  shows that
\begin{equation}\label{newGamma}
\Gamma(t,h) = \fv(h)+ d \fA(h) \widetilde \gamma(t+\overline h, h)
\end{equation}
where
$\fv(h)$ is a vector in $\bbR^d$ with coordinates
$\fv_k(h)=\sum_{\nu=1}^d (\ka_\nu(h)-\overline h)^k$ and $\fv_d(h)=0$,
and $\fA(h)$ is a $d\times d$ matrix with
\begin{equation*}
\fA_{ij}(h)=\begin{cases}
1, \quad i=j,
\\
0,  \quad i>j
\\
d^{-1}
\sum_{\nu=1}^d \tfrac{(\ka_\nu(h)-\overline h)^{j-i}}{(j-i)!}, \quad
i<j\le d-1\\0,\quad  i<d, j=d.
\end{cases}
\end{equation*}
 Finally $\widetilde \gamma(s,h)=
(s,\dots, \tfrac{s^{d-1}}{(d-1)!}, \widetilde \phi(s; h))$ with
$$\widetilde \phi(s;h)= \frac 1d\sum_{i=1}^d \phi(s-\overline
h+\kappa_i(h)).$$
The function $\widetilde \phi$ and the curve
$\widetilde \gamma(t,h)$ are defined on
$[a(h), b(h)]\subset [0,1]$
where $a(h)=a+\overline h$ and
$b(h)=b-\ka_d(h)+\overline h$.
It is now crucial to note that for $\phi\in \cK_{a,b,M}(\sigma)$ and fixed
$h$ the offspring function $\widetilde \phi \equiv \widetilde \phi(\cdot;h)$
belongs to
$\cK_{a(h),b(h),M}(\sigma/d)$. This follows  from
\eqref{twoJs}, \eqref{integraldarstellung} for the function $\widetilde \phi$.
Indeed
the nonnegativity of $\Psi _d$
imply that if
$\widetilde h\in (\mathbb R _+ )^{d-1}$ satisfies $\ka_d(\widetilde h)\le
b(h)-a(h)$, then
\begin{align*}
J_{\widetilde \phi(\cdot ,h)} (t, \widetilde h)&=
\int\nolimits_{t}^{t+\ka_d(\widetilde h)}
\Psi_d (u-t;\tilde h) \frac{1}{d} \sum_{i=1}^{d}
\phi ^{(d)}
(u-\overline h +\ka_i(h)) \, du
\\&\ge
\frac{\sigma}{d}\, v(h)\, \sum_{i=1}^d \Bigl(\prod_{j=1}^d
\phi^{(d)}(t-\overline h +\ka _i (h)+\ka _j (\tilde{h}))
\Bigr)^{1/d}
\\&\ge
\frac{\sigma}{d}\, v(h)\, \Bigl(\prod_{j=1}^d
\phi^{(d)}(t-\overline h +\ka _d (h)+\ka _j (\tilde{h}))
\Bigr)^{1/d}
\\&\ge
\frac{\sigma}{d}\, v(h)\,  \prod_{j=1}^d
\Bigl(\frac{1}{d}\sum_{i=1}^d
\phi^{(d)}\bigl( t-\overline h +\ka _i (h)+\ka _j (\tilde{h})\bigr)
\Bigr)^{1/d}.
\end{align*}
Here the first inequality follows from \eqref{integraldarstellung} and
$\phi\in \cK_{a,b}(\sigma )$. The last  inequality  shows that
$\tilde{\phi}(\cdot
;h)\in \cK_{a(h),b(h)}(\sigma )$; it  follows from the fact that
$\phi^{(d)}$
is nondecreasing.


Now let $g_h$ be defined  by
$\widehat g_h(\xi)=\widehat g(\fv(h)+d \fA(h) \xi)$. Then  because of the
 unimodularity of $\fA(h)$ we have
$\|g_h\|_{q'}=d^{d/q'}\|g\|_{q'}$. Also if $g$ is supported in
$B_R$ then $g_h$ is supported in the ball of radius $Rd^3$
(observe that all the entries of $\fA(h)$ are at most $d$).

Comparing a geometric to an arithmetic mean we see that
\begin{align*}
&\Big(\int_a^{b-\ka_d(h)} \big|
\widehat g(\Gamma(t,h)) \big|^{p'} H(t,h)^{1/d} dt\Big)^{1/p'}
\\&\le
\Big(\int_a^{b-\ka_d(h)} \big|
\widehat g(\Gamma(t,h)) \big|^{p'}
\Big(\frac{1}{d}\sum_{i=1}^d \phi^{(d)}(t+\ka_i(h))\Big)^{2/(d^2+d)}
dt\Big)^{1/p'}
\\&=
\Big(\int_{a+\overline h}^{b-\ka_d(h)+\overline h} \big|
\widehat g_h(\widetilde \gamma(s,h)) \big|^{p'}
\big(\widetilde \phi^{(d)}(s;h)
\big)^{2/(d^2+d)}
ds\Big )^{1/p'}
\\&\le
\frac{\fc/d}{\sigma/d}\, \cA(R d^3, M, \sigma/d) \|g_h\|_{q'}
=\frac{\fc}{\sigma} d^{d/q'} \cA(R d^3, M, \sigma/d) \|g\|_{q'} .
\end{align*}
By duality this also implies \eqref{uniform2}.
\end{proof}

We now proceed exactly as in the proof of Proposition 6.1 in \cite{BOS}.
We first have, by an application of
Plancherel's theorem and the  change of variable
$(t,h) \mapsto
\Gamma(t,h)$
\begin{equation} \label{planchappl}
\Big\|\int S_{R,h}[F]dh\Big\|_2
\le C \Big( \iint \big|F(t,h)H(t,h)J(t,h)^{-1/2}\big|^2 dt
\, dh \Big)^{1/2};
\end{equation}
the change of variable  can be justified as in
\cite{DM1}, p. 549.

 Replacing $F$ with $FH^{(d-1)/d}$ in \eqref{uniform2} and then integrating
with respect to $h$ now yields, according to Minkowski's inequality, the
estimate
\begin{multline} \label{interp0}
\Big\|\int S_{h}[F]dh\Big\|_{L^q (B_R )}
\\
\le
d^{d/q'}\,\fc\, \sigma^{-1} \,
\cA(R d^3, M,  \fc /d)
 \int\Big( \int |F(t,h) H(t,h)^{1-\frac{1}{d}+\frac{1}{dp}}|^p
dt \Big)^{1/p}dh .
\end{multline}

By analytic interpolation of \eqref{interp0}
and \eqref{planchappl} one  obtains
%
\begin{multline}\label{ineq}
\Big\| \int S_{R,h}[F] dh\Big\|_{L^s (B_R )}
\le C  \Big(\frac{\fc}{\sigma}\cA(R d^3, M, \sigma/d)\Big)^{1-\vth}
\\
\times\Big(\int \Big(\int
\big|F(t,h)H(t,h)^{ \eta } J(t, h )^{-\vartheta/2 }\big|^{B(\vth)}
dt\Big)^{A(\vth)/B(\vth)} dh \Big)^{1/A(\vth)}
\end{multline}
where $0\le \vartheta \le 1$ and $A,B,s,\eta$ are defined by
\begin{equation}\label{vth}
\begin{aligned}
&\frac 1{A(\vth)}=  1-\frac \vartheta {2}, \quad
&&\frac{1}{B(\vth)}= \frac{1}{p}+\vth(\frac{1}{2}-\frac {1}{p}),
\\
&\frac{1}{s(\vth)}=
\frac{1-\vartheta}{q}+\frac{\vartheta}{2}, \quad
&&\eta(\vth) =1- \frac{d+1}{2q}(1-\vth).
\end{aligned}
\end{equation}

Now let
\begin{equation}\label{vthp}
 \vth(p)= \frac{4(d-1)}{(d+1)d p'-4} = \frac{2(d-1)}{q-2}
\end{equation}
and let $A_p=A(\vth(p))$, $B_p=B(\vth(p))$, $s_p=s(\vth(p))$ and
$\eta=\eta(\vth(p))$. Then
$\eta_p -(d+1)\vth/4=  1/p$ and
$s_p = q/d= (d+1)p'/2.$
As $\phi\in \cK_{a,b,M}(\sigma)$ we may use the crucial inequality
$J_\phi(t)\ge \sigma v(h) H^{(d+1)/2}(t,h)$
and obtain
\begin{multline} \label{vthpinequality}
\Big\|\int S_{h}[F] dh\Big\|_{L^{q/d}(B_R)}
 \le C \sigma^{-\vth(p)/2}
( \fc \,\sigma^{-1}\cA(R d^3, M, \sigma/d))^{1-\vartheta(p) } \times \\
\Big(\int \Big(\int
\prod\limits_{j=1}^d \big| F(t,h) H(t,h)^{\eta_p- \frac{d+1}{4}\vth(p)}
 \big|^{B_p} dt\Big)^{A_p/B_p}
v(h)^{1-A_p}dh \Big)^{1/A_p}.
\end{multline}

We are now in the position to apply an inequality by Drury and Marshall
\cite{DM1}
for multilinear operators involving Vandermonde's determinant,
see also \cite{BOS} for an exposition.
To state this
let
$$\frak V[f_1,\dots, f_d](t,h)
:= v(h)^{-1} \prod_{i=1}^d f_i(t+\kappa_i(h))
$$
and  $L_v^A(L^B)$
denote the weighted mixed norm space
consisting of func\-tions $(t,h)\mapsto G(t,h)$ with
$\|G\|_{L_v^A(L^B)}=(\int \|G(\cdot,h)\|_B^A
v(h)dh)^{1/A}<\infty$. One assumes that   $1<A<\frac{d+2}d$,
$1<A\le B<\frac{2A}{d+2-dA}$, and sets $\sigma =2/(d+2-dA)$. For
$l=1,\dots, d$ let $Q_l$ denote the  point in $\bbR^d_+$ for which
the $j^{\text{th}}$ coordinate is  $(\sigma A)^{-1}$, if $j\neq l$
and the $l^{\text {th}}$ coordinate is $B^{-1}$. Let $\Sigma(A,B)$
be the $d-1$ dimensional closed convex hull of the points
$Q_1,\dots, Q_d$. Then the inequality
\begin{equation}\label{vandineq}
\big\|\frak V[f_1,\dots, f_d]\big\|_{L_v^A(L^B)} \le
C\prod_{i=1}^d \|f_i\|_{L^{p_i,1}}
\end{equation}
holds for all  $(p_1^{-1}, \dots, p_d^{-1})\in \Sigma(A,B)$.

We apply this inequality to the right hand side of
 \eqref{vthpinequality} to obtain
 \begin{multline*}
\Big\|\int S_h [F]dh\Big\|_{L^{q/d}(B_R )}\\
\le
C(d,p)\, \sigma ^{-\vth (p)/2}\Bigl(\fc\,\sigma^{-1}\,
\cA(R d^3,M,
\fc /d)\Bigr) ^{1-\vartheta(p) }
\prod_{j=1}^d \|f_j w^{1/p}\|_{L^{p_j ,1}}
\end{multline*}
whenever $(p_1^{-1},\dots ,p_d^{-1})\in \Sigma (A_p ,B_p )$.
Summing over the permutations $\pi\in \fS_d$ then yields
\begin{multline}
\label{prodineq}
\Big\|\prod_{i=1}^d Tf_i \Big\|_{L^{q/d}(B_R )}\le\\
C(d,p)\, \sigma ^{-\vth (p)/2}\Bigl(\fc\,\sigma^{-1}\, \cA(R d^3,
M,\fc /d)\Bigr) ^{1-\vartheta(p) }
\prod_{j=1}^d \|f_j w^{1/p}\|_{L^{p_j ,1}}.
\end{multline}

We now use applications of H\"older's inequality and
 Christ's multilinear trick for the $q_d$-linear expression $\prod_{i=1}^{q_d} T f_i$,
exactly as in  \S6 of \cite{BOS}.
This  yields
\begin{multline*}
\Big\|
\prod_{i=1}^{q_d} T f_i\Big\|_{L^{q/q_d}(B_R )}
\\
\lc \sigma^{-q_d\vth(p)/2d}
( \fc\,\sigma^{-1}\, \cA(R d^3, M, \sigma/d))^{(1-\vartheta(p))q_d/d}
\prod_{i=1}^{q_d}\big\|f_i |\phi^{(d)}|^{\frac{2}{(d^2+d)p}}\big\|_{L^{p,q_d}}.
\end{multline*}
Since $p<q_d<q$ this implies (for $f_i\equiv f$)
\begin{multline}
\label{outcome}
\|T f\|_{L^q (B_R )} \le \\C(d,p,q) \sigma^{-\vth(p)/2d}
( \fc\,\sigma^{-1}\, \cA(R d^3, M, \sigma/d))^{(1-\vartheta(p))/d }
\big\|f |\phi^{(d)}|^{\frac{2}{(d^2+d)p}}\big\|_p
\end{multline}
provided that $\phi\in \cK_{a,b,M}$ for some $M<\infty$.
Observe that from the definition of $\cA$ we get
 $$\cA(R d^3, M, \fc/d))\le C_{d,p} \cA(R, M, \fc)$$ and thus
\eqref{outcome} implies
$$\cA(R,M,\fc)\le C(d,p)
\cA(R,M,\fc)^{(1-\vth(p))/d} \sigma^{-\vth(p)/2d}
$$
which by \eqref{vthp} yields \eqref{cAclaim}. \qed

\section{Examples of curves covered by Theorem \ref{flatthm}}\label{examples}

\subsection{}
 Condition \eqref{techcond} (and {\it a fortiori} condition \eqref{techcondnew}) holds
 for $\phi(t)=t^\beta$ and the required
monotonicity of the first $d$ derivatives holds if $\beta>d-1$.

\subsection{} Consider the function
$\phi(t)= \exp(-t^{-\beta})$ for $t>0$.
Then
induction  shows that
$\phi^{(d)}(t) =\beta^d e^{-t^{-\beta}}
t^{- d(\beta+1)} \big(1+\sum_{j=1}^d a_{j,d} t^{j\beta}\big)
$
and the  coefficients satisfies the recursive relation
$a_{k,d+1}=\beta^{-1}  a_{k,d} -
a_{k-1,d}(d+1-k+d/\beta)$ if $k
\le d-1$ and
$a_{d,d+1}= -
a_{d-1,d}(1+d/\beta)$ if $k=d$.
It is obvious that if $A>1$, then condition \eqref{techcond} is satisfied on
a (small) interval $(0,c(A))$.

\subsection{} \label{flex}
Suppose that $\big(\prod_{j=1}^d g(s_j)\big)^{1/d} \le  g\big(\root \uproot 1d \of{s_1\cdots s_d}\big)$
for  $0<s_1\le s_2\le \dots\le s_d<\infty$, and $g$ is nonnegative and increasing.
Set $f_g(s)= \exp(-1/g(s))$.  Then we also have for $\overline s=(\prod_{i=1}^d s_j)^{1/d}$
\begin{align*}
f_g(\bar s) &=\exp \big(-1/g(\overline s )\big)\geq
\exp \big(-   ( \prod_{j=1}^d 1/g(s_j ))^{1/d}\big)
\\
&\ge \exp \Bigl(-\frac{1}{d}\sum_{j=1}^d \frac{1}{g (s_j )}
\Bigr)=\Bigl(\prod\limits_{j=1}^d f_g(s_j )\Bigr)^{1/d}.
\end{align*}
Thus if the first $d$ derivatives of a function $\phi$ are nonnegative
and increasing on $(0,\infty)$ and if $\phi$ satisfies \eqref{techcond}
with $A=1$
then the same  conditions are satisfied by
$\psi(t)=\int_0^t (t-u)^{d-1} \exp(-1/\phi^{(d)}(u)) du$.
As mentioned in the introduction this leads to a sequence of progressively
flatter functions mentioned following the statement of Theorem \ref{flatthm}.

\subsection{}\label{logphid}
 Similarly, suppose that $\big(\prod_{j=1}^d g(s_j)\big)^{1/d} \le
g\big(\root\uproot 1d \of{s_1\cdots s_d}\big)$
for $0\le a<s_1\le s_2\le \dots\le s_d<b$. Assume also
that $g(s)>e$ if $s\in (a,b)$. Then
\begin{multline*}
\Big(\prod_{j=1}^d \log (g(s_j))\Big)^{1/d}
\le \frac{1}{d}\sum_{j=1}^d \log (g(s_j ))
\\
=\log\Big(\prod_{j=1}^d
g(s_j)\Big)^{1/d}\le \log (g(\root\uproot 1d \of{s_1\cdots s_d})).
\end{multline*}
Again if
$\psi(t)=\int_a^t (t-u)^{d-1} \log(\phi^{(d)}(u)) du$,
if $\phi^{(d)}(s)>e$ and $\phi^{(d)}$ is  non\-decreasing on $(a,b)$
then condition
\eqref{techcond} with $A=1$ for $\phi$  implies
\eqref{techcond} with $A=1$
for $\psi$.

\section{Proof of Theorem \ref{euclthm}}\label{eucl}

First assume that \eqref{hypothesis2} holds. We will establish
\eqref{conclusion2}. For $\la>1$ define
$$
T_{\la} f(x) = \chi(x) \int_a^b f(t) e^{-i\la \inn{x}{\gamma(t)}} dt,
$$
where $\chi$ is  the characteristic function of a set of diameter $1$.

\medspace

\noindent{\bf Definition.}
\textit{ For $-\infty <a<b<\infty$ and $\sigma>0$,
let $\cC_{a,b}(\sigma)$ be the
class  of all real-valued functions $\phi$ defined on $(a,b)$ for which}

\textit{(i) $\phi\in C^d((a,b))$ and the derivatives
$\phi',\dots, \phi^{(d)}$ are nonnegative and nondecreasing
on $(a,b)$, and }

\textit{(ii) the inequality
\begin{equation}\label{phicond}
\phi ^{(d-1)}(s)-\phi ^{(d-1)}(t)\ge \sigma ^{-\frac{1}{\alpha}}
(s-t)^{\frac{1}{\alpha}+1-\frac{d(d+1)}{2}}
\end{equation}
holds for all $s$ and $t$ such that $a< t< s< b$.
}

\medskip

\noindent With $q=1+1/\alpha$ and for $\lambda >1$, $\sigma >0$ and
large $r$, define
\begin{equation*}
\cB\equiv \cB(\la ,\sigma,r) :=
\la^{d/q}\sup_{\substack{
\phi\in \cC_{a,b}(\sigma)\\-r\le a<b \le r
} }
\sup_{ \|f\|_{L^{q}((a,b))}\le 1}
\|T_{\la}f\|_{L^{q,\infty}(\mathbb R ^d )}.
\end{equation*}

\noindent By duality and Lemma \ref{Lemma1} below, \eqref{conclusion2}
is a consequence of the following estimate
\begin{equation}\label{conclusion}
\cB(\la ,\sigma,r)\le C(d,\alpha ) \, \sigma ^{\frac{1}{1+\alpha}}.
\end{equation}

\begin{lemma}\label{Lemma1}If \eqref{hypothesis2} holds for all
parallelepipeds $E$ in $\mathbb R^d$ then the inequality
\begin{equation*}
B^{-\frac{1}{\alpha}}\,
(s-t)^{\tfrac{1}{\alpha}+1-\frac{d(d+1)}{2}}\le
\phi ^{(d-1)}(s)-\phi ^{(d-1)}(t)
\end{equation*}
holds whenever $a<t<s<b$.
\end{lemma}
We shall give the proof in \S\ref{lemmas}.

To begin the proof of \eqref{conclusion}, fix $a$, $b$, $\sigma$,
and $\phi\in\cC_{a,b}(\sigma )$ and then define
\begin{align*}
M_{\la}&(f_1, \cdots, f_d)(x)=\prod_{j=1}^d T_{\la}f_j (x)\\&=
\chi(x)\int_{\mathbb R^{d-1}}\int_{\cI_h}e^{-i\la
\inn{x}{\sum_{j=1}^d \gamma(s+h_j )} }\prod_{j=1}^d f_j (s+h_j)\, ds
\, dh_1 \cdots dh_{d-1} ,
\end{align*}
where our convention now is that $h_d =0$ and $\cI_h$ is the
(possibly empty)
intersection of the $d$ intervals $(a-h_j ,b-h_j )$. In what follows we
will further simplify the notation by writing $h=(h_1 ,\dots ,h_{d-1} )$
and
$\Gamma (s,h)=\sum_{j=1}^d \gamma(s+h_j )$. With an eye to decomposing the
multilinear operator $M_{\la}$ we define
$$u(h)
=\prod_{1\le i<j\le d}|h_i -h_j |
=h_1\cdots h_{d-1}\prod_{1\le i<j\le d-1}|h_i -h_j |
$$
and
$$K(h)=u(h)\Bigl(\sup_{1\le i<j\le d}|h_i -h_j
|\Bigr)^{\frac{1}{\alpha}-\frac{d(d+1)}{2}}.$$
Note that  $K$ is homogeneous of degree $\alpha^{-1}-d.$
Now, for $m \in \bbZ$, let $$S_m =\{h\in \bbR^{d-1}:2^{-m-1}<K(h)\le 2^{-m}\}$$
and, following \cite{baklee},  define
\begin{equation*}\label{}
M_{\la, m} (f_1, \cdots, f_d)(x) = \chi(x) \int_{S_m} \int_{\cI_h}
e^{ -i\la\inn{x}{\Gamma (s,h)}
} \prod_{j=1}^d f_j(s+h_j)\, ds \, dh .
\end{equation*}
We will need to observe that
\begin{equation}\label{S_m}
m_{d-1}(S_m )\le C(d,\alpha )\, 2^{-m(d-1)\alpha /(1-d\alpha )}.
\end{equation}
By homogeneity, it is enough to check that $m_{d-1}(\{h:K(h)\le 1\})\le
C(d)$. Since $\alpha\le 2/(d^2 +d)$,
$$\{h:K(h)\le 1\}\subset (\{h:u(h)\le 1\}
\cup\{h:\sup |h_i |\le 1\},$$
and so it is enough to check that
\begin{equation}\label{uset}
m_{d-1}(\{h:u(h)\le 1\} )\le C(d).
\end{equation}
But it follows from  \cite{DM1} (see (i) of Proposition 2.4 in \cite{BOS}) that
$$
m_{d-1}(\{ h:
0\le h_1 \le\cdots\le h_{d-1}\,;\ u(h)\le 1 \})
\le C(d)
$$
and so
$ m_{d-1}(\{ h:0\le h_j \, ;\ u(h)\le 1 \})
\le C'(d).$ Since
$$
\prod_{1\le i\le j\le d}\big| |h_i |-|h_j |\big |
\le \prod_{1\le i\le j\le d}|h_i -h_j |=u(h),
$$
\eqref{uset} follows.

Now considerations similar to those which lead to \eqref{newGamma} show
that
$$
\Gamma(s,h) = \fv(h)+ d \fA(h) \widetilde \gamma(s+\overline h, h)
$$
where $\fv(h)$ is a vector, where $\fA(h)$ is a matrix with
entries $1$ on the diagonal and $0$ below, where $\overline h
=\sum_{j=1}^d  h_j
/d$, and where
$$\widetilde{\gamma}(s,h)=\Bigl(s,\frac {s^2}{2},\dots
,\frac{s^{d-1}}{(d-1)!},\widetilde{\phi}(s,h)\Bigr)$$ with
$$
\widetilde{\phi}(s,h)=\frac{1}{d}\sum_{i=1}^d \phi (s-\overline h
+h_i ).
$$
Since \eqref{phicond} holds for $\phi$, it holds as well for each
$\widetilde{\phi}(\cdot ,h)$. Therefore we have the estimate
$$
\Big\| \chi\,\int_{\cI_h}e^{-i\la \inn{\cdot}{\Gamma (s,h)}}f(s)\, ds
\Big\|_{L^{q,\infty}(\mathbb R ^d)}\le \la ^{-d/q}\cB (\lambda
,\sigma,r)\, \|f\|_{L^q (\cI_h )}.
$$
Taking  \eqref{S_m} into consideration,
an application of Minkowski's
inequality thus yields
\begin{multline}\label{Eq}
\| M_{\la, m} (f_1, \cdots, f_d)\|_{L^{q,\infty}(\mathbb R ^d )}
\\\le
C(d,\alpha ) \la^{-d/q}
\cB (\lambda ,\sigma,r ) \,  2^{-m [(d-1)\alpha/(1 - d \alpha )]} \| f_d
\|_q
\prod_{j=1}^{d-1} \| f_j\|_{\infty} ,
\end{multline}
where $\|\cdot\|_q$ stands for the norm in $L^q (a,b)$.

Let $J(s,h)$ stand for the absolute value of the Jacobi-determinant of
the transformation $(s,h) \mapsto
\Gamma(s,h)$ (defined on $\{(s,h):s\in \cI_h \}$).
To obtain an $L^2$ estimate for $M_{\lambda ,m}$ we will need the
following inequality:
\begin{equation}\label{Jacobianest}
J(s,h) \ge c(d)\, \sigma ^{-1/\alpha} K(h).
\end{equation}
This inequality follows from \eqref{phicond}
and the next lemma
whose proof is given in
\S\ref{lemmas}.

\begin{lemma}\label{Lemma2} Suppose the inequality \begin{equation}
\label{phiest}
 c \,(s-t)^{\rho}\le\phi ^{(d-1)}(s)-\phi ^{(d-1)}(t)
\end{equation} for some $\rho >0$ and
for $a<t<s<b$. Then there is also the inequality
\begin{equation*}
c(d) \,c \, u(h)\Bigl(\sup_{1\le i<j\le d}|h_i -h_j
|\Bigr)^{\rho -1}\le J(s,h )
\end{equation*}
whenever $s\in \cI_h$.
\end{lemma}

Now the transformation $(s,h)\mapsto \Gamma (s,h)$ is at most $d!$ to one
a.e.,
so
$$
\|M_{\lambda ,m}(f_1 ,\dots ,f_d)\|_{L^2 (\mathbb R ^d)}^2
\le d!\, \int_{S_m}\int_{\cI_h} \Big| \prod_{j=1}^d f_j (s+h_j )\Big| ^2
\frac{1}{J(s,h)}ds\, dh .
$$
Applying \eqref{Jacobianest} and recalling \eqref{S_m},
we obtain
\begin{multline}\label{E2}
\| M_{\la, m} (f_1, \cdots, f_d)\|_{L^2 (\mathbb R ^d )} \\\le C(d,\alpha
)\,
\la^{-d/2}\, \sigma^{1/2\alpha}\, 2^{m
[1-(2d-1)\alpha]/[2(1-d\alpha)]} \| f_d \|_2 \prod_{j=1}^{d-1} \|
f_j\|_{\infty} .
\end{multline}

Interpolating the estimates \eqref{Eq} and \eqref{E2} yields that
\begin{multline}\label{Mlamest}
\| M_{\la,m} (f_1, \cdots, f_d)\|_{L^{q/d,\infty}(\mathbb R ^d )}
\\ \le C(d,\alpha )\,
\la^{-d^2/q} \, \sigma^{(d-1)/(1-\alpha )}\,
\cB (\lambda ,\sigma,r)^{\delta(\alpha)}
\| f_d\|_{q/d}
\prod_{j=1}^{d-1}
\| f_j\|_{\infty},
\end{multline}
with $$\delta(\alpha)  = \frac{1-(2d-1)\alpha}{1-\alpha} \in (0,1).$$
If one  uses Bourgain's interpolation argument in \cite{Bo1}
(see also the appendix of \cite{CSWW}) then one  actually
obtains an estimate for the sum $M_\la=\sum_m M_{\la,m}$, namely,
\begin{multline}\label{Mlaest}
\| M_{\la} (f_1, \cdots, f_d)\|_{L^{q/d,\infty}(\mathbb R ^d )}
\\ \le C(d,\alpha )\,
\la^{-d^2/q} \, \sigma^{(d-1)/(1-\alpha )}\,
\cB (\lambda ,\sigma,r)^{\delta(\alpha)}
\| f_d\|_{q/d,1}
\prod_{j=1}^{d-1}
\| f_j\|_{\infty}.
\end{multline}
To arrive at \eqref{Mlaest} it suffices to prove this bound for
$f_d=\chi_U$, the characteristic function of a measurable set $U$.
One then uses \eqref{E2} to estimate the size of the set where
$|\sum_{2^m\le \beta} M_{\la,m} (f_1,\dots, f_{d-1},\chi_U)|\ge
s$, and one uses  \eqref{Eq} to estimate the size of the set where
$|\sum_{2^m> \beta} M_{\la,m} (f_1,\dots, f_{d-1},\chi_U)|\ge s$;
here $\beta>0$ will be suitably chosen. This leads to
\begin{align*}
&m_d\big(\big\{x: \big|\sum_m M_{\la,m} f(x)\big| >2s\big\}
\big)
\\
&\le\,\,\la^{-d}|U|\Big[
C(d,\alpha)^q s^{-q} \cB(\la, \sigma,r)^q \beta^{-\frac{(d-1)\alpha q}{1-d\alpha}}
\prod_{i=1}^{d-1}\|f_i\|_\infty^q
\\
&\qquad\qquad+
C(d,\alpha)^2 s^{-2} \sigma^{1/\alpha} \beta^{-\frac{(1-(2d-1)\alpha}{1-d\alpha}}
\prod_{i=1}^{d-1}\|f_i\|_\infty^2\Big],
\end{align*}
and the estimate
\eqref{Mlaest} follows by choosing the optimal $\beta$.
\eqref{Mlaest} gives
\begin{multline}\label{Est}
\Big\| \prod_{j=1}^d T_{ \la} f_j
\Big\|_{L^{\frac{q}{d},\infty}(\mathbb R ^d )} \\
\le C(d,\alpha )\,
\la^{-d^2/q}
\, \sigma^{(d-1)/(1-\alpha )}\,
\cB (\lambda ,\sigma,r )^{\delta(\alpha)} \, \|
f_1\|_{\frac{q}{d},1} \prod_{j=2}^d \| f_j\|_{\infty},
\end{multline}
and if  we take for all $f_j$ the same characteristic function of a set
we  also get
\begin{equation}\label{Est2}
\left\| T_{ \la} f \right\|_{q,\infty} \le C(d,\alpha )\,  \la^{-d/q}
\, \sigma^{(d-1)/(d-d\alpha )}\,
\cB(\lambda ,\sigma,r )^{\delta(\alpha)/d} \| f\|_{q ,1}.
\end{equation}

Now fix an integer $N > q$.
Applying a version of H\"older's inequality
(see (2.1) in \cite{BOS}) and permuting the functions,
\eqref{Est} and \eqref{Est2} yield
\begin{multline*}\label{}
\Big\| \prod_{j=1}^N T_{\la} f_j \Big\|_{L^{q/N,\infty}(\mathbb R
^d )}
\\
\le C(d,\alpha )\,  \la^{-N d/q} \, \sigma^{N(d-1)/(d-d\alpha )}\, \cB
(\lambda ,\sigma,r )^{N \delta(\alpha)/d} \prod_{j=1}^N \|
f_j
\|_{L^{q_j, 1}}
\end{multline*}
when $(q_1^{-1}, \cdots, q_N^{-1})$ is one of the $N$ points $Q_j$
in $\bbR^N$ defined as follows: $Q_1$ is the point with the first
component
$d/q$, the next $d-1$ components $0$, and the remaining $N-d$
components equal to $1/q$; $Q_2$ is obtained by shifting the
components of $Q_1$ to the right by one and moving the last
component to the front; etc. Here $L^{\infty, 1}$
should be interpreted as $L^{\infty}$.
Applying Christ's  multilinear trick (for multilinear operators
with values in the   quasi-normed $q/N$-convex space $L^{q/N,\infty}$,
 see
Proposition 2.3  in \cite{BOS} and also \cite{janson}), these estimates yield
\begin{multline*}\label{}
\Big\| \prod_{j=1}^N T_{ \la} f_j \Big\|_{L^{q/N,\infty}(\mathbb R ^d )}
\\
\le C(d,\alpha )\,
 \la^{-N d/q}
\, \sigma^{N(d-1)/(d-d\alpha )}\,
\cB (\lambda ,\sigma,r )^{N \delta(\alpha)/d}
\prod_{j=1}^N \|
f_j
\|_{L^{q_j, r_j}}
\end{multline*}
when $(q_1^{-1}, \cdots, q_N^{-1})$ is in the interior of the
convex hull $\Sigma$ of $Q_1, \cdots, Q_N$ and when the $r_j \in [1,
\infty]$ satisfy $\sum_{j=1}^N 1/r_j = N/q$. Note that the point
$(1/q, \cdots, 1/q)$ is the center of $\Sigma$.
Hence, taking $f_j = f$ and $q_j = r_j =q$, we obtain
\begin{equation*}\label{}
\| T_{ \la} f\|_{L^{q,\infty}(\mathbb R ^d )} \le C(d,\alpha )\,
\la^{-d/q}
\, \sigma^{(d-1)/(d-d\alpha )}\,
\cB (\lambda ,\sigma,r )^{\delta(\alpha)/d} \| f\|_{L^{q}} .
\end{equation*}
Therefore, by the definition of $\cB (\lambda,\sigma,r )$, we have
\begin{equation*}\label{}
\cB (\lambda ,\sigma,r ) \le C(d,\alpha )\, \sigma^{(d-1)/(d-d\alpha )}\,
\cB
(\lambda ,\sigma,r )^{\delta(\alpha)/d}.
\end{equation*}
Recalling the definition of $\delta$, some algebra yields
\eqref{conclusion}. Thus \eqref{conclusion2} is established.

Now for the converse, we  assume that \eqref{hypothesis3} holds
with $1/P'=\alpha /Q$ and will show
that \eqref{hypothesis2} holds with $B$ replaced by $C(d,p)\, B$.
Fix an $f\in C_c^\infty (\mathbb R ^d )$ with $f$ nonnegative and equal to
$1$ on $[0,1]^d$. Consider a parallelepiped
$$
E=x_0 +\big\{\sum_{j=1}^d t_j x_j :0\le t_j \le 1\big\}
$$
and fix a linear isomorphism $T$ of $\mathbb R ^d$ which satisfies
$$
T([0,1]^d)=\big\{\sum_{j=1}^d t_j x_j :0\le t_j \le 1\big\}.
$$
Let $g$ be defined by $\widehat g (x)=f(T^{-1}(x-x_0 ))$ so that
$\widehat g$ is nonnegative and equal to $1$ on $E$ . Then
a computation shows $\|g\|_{L^{P}(\mathbb R ^d )}=
m_d (E)^{1/P'}
\ \|\widehat f\|_{L^{P}(\mathbb R ^d )}$. If \eqref{hypothesis3} holds
then it follows that
$$
\lambda (E)^{1/Q}\le \Big(\int_a^b
 |\widehat g (\gamma(t))|^{Q}  dt\Big)^{1/Q} \le
B^{\frac{1}{Q} } m_d (E)^{1/P'}
\ \|\widehat f\|_{L^{P}(\mathbb R ^d )}.
$$
Since $1/P'=\alpha /Q$ this yields \eqref{hypothesis2}
with $B$ replaced by $\|\widehat f \|_{L^{P}(\mathbb R ^d )}^{Q}\, B$
and therefore completes the proof of Theorem \ref{euclthm}.

\section{Proofs of Lemma \ref{Lemma1} and Lemma \ref{Lemma2}}\label{lemmas}

\begin{proof}[{\bf Proof of Lemma \ref{Lemma1}}]
Write $s=t+h$ and let $E_{d-2}$ be the parallelogram in
$\mathbb R ^2$ with vertices
\begin{align*}
P_1&= (t, \phi^{(d-2)}(t)),  &&P_2=P_1-\rho  e_2
\\
P_3&= (t+h, \phi^{(d-2)}(t+h)),  &&P_4=P_3+\rho  e_2
\end{align*}
where
$\rho= h\phi^{(d-1)}(t+h)+ \phi^{(d-2)}(t)-\phi^{(d-2)}(t+h)\ge 0$, so that
$\phi^{(d-1)}(t+h)$ is the slope of the line segments
$\overline{P_2P_3}$ and  $\overline{P_1P_4}$.
Then (as a sketch will show)
\begin{equation}\begin{aligned}\label{k=2calculation}
m_2 (E_{d-2})&\leq 2\int\nolimits_t^{t+h}\bigl(\phi ^{(d-2)}(t)+\phi
^{(d-1)}(t+h)(s-
t) -\phi ^{(d-2)}(s)\bigr)\, ds
\\&=
2\int\nolimits_t^{t+h}\int\nolimits_t^{s}\bigl(\phi ^{(d-1)}
(t+h)-\phi ^{(d-1)}(u)\bigr)\, du\, ds
\\ &\leq
2\int_t^{t+h}\int\nolimits_t^{s}\bigl(\phi ^{(d-1)}(t+h)-\phi
^{(d-1)}(t)\bigr)\, du\, ds
\\&=
h^2 \bigl(\phi ^{(d-1)}(t+h)-\phi ^{(d-1)}(t)\bigr).
\end{aligned}
\end{equation}
We now prove the following

{\it Claim:} For $2\le k\le d$,
\begin{equation}
\label{claimestimate}
\{\gamma ^{(d-k)}(s):t\leq s\leq t+h\}\subset
\begin{cases}
\{e_{d-k}\}\times E_{d-k}, \quad &2\le k\le d-1,
\\
E_0, &k=d,
\end{cases}
\end{equation}
where $\{e_1,\dots, e_d\}$ is the standard basis in $\Bbb R^d$ and
$E_{d-k}$ is a parallelepiped in $\bbR^k$ with
\begin{equation}
\label{edkest}
m_k(E_{d-k}) \le h^{\frac{k^2+k-2}{2}}
\big(
\phi^{(d-1)}(t+h)-
\phi^{(d-1)}(t)\big).
\end{equation}

The above calculation \eqref{k=2calculation}
verifies this claim for $k=2$, and all $d\ge 2$.
We argue by induction on $k$ and assume $3\le k\le d$ and that the induction hypothesis is true for $k-1$.

Now suppose $s\in [t,t+h]$. Then
$$
\gamma ^{(d-k)}(s)-\gamma ^{(d-k)}(t)=\int_t^s \gamma ^{(d-k+1)}(u)
\, du
$$
belongs to
\begin{align*}&O_{d-k}\times (s-t)
\big(\{e_{d-k+1}\}\times E_{d-k+1}\big)
\\&\subset
O_{d-k}\times \{u(1,x)\in \bbR\times  \bbR^{k-1}: 0\le u\le h, x\in
E_{d-k+1}\}
\end{align*}
where $O_{d-k}$ denotes the origin in $\bbR^{d-k}$ and where $O_{d-k}$ is omitted  if $ k=d$.

Let $ x_0$ be any point of $E_{d-k+1}$ in $\mathbb R ^{k-1}$. The set
$$\widetilde E_{d-k}\defeq\{(1,x)-v(1, x_0 ):x\in E_{d-k+1},\ 0\le v\le 1\}$$
is  a
parallelepiped in $\mathbb R ^k$ which satisfies $m_{k}(\widetilde E_{d-k})=m_{k-1}
(E_{d-k+1})$, which contains $O_k$ and
$\{(1,x):x\in E_{d-k+1}\}$, and which therefore (by convexity) contains
$$
\{u(1,x):0\le u\le 1,\ x\in E_{d-k+1}\}.
$$

Thus, with
$$
E_{d-k}\defeq \big\{(t,\dots, \tfrac{t^{k-1}}{(k-1)!}, \phi ^{(d-k)}(t))+uy:
0\le u\le h,\ y\in \widetilde E_{d-k}\},
$$
we have
$$
\{\gamma ^{(d-k)}(s):t\leq s\leq t+h\}\subset
\begin{cases} \{e_{d-k}\}\times E_{d-k},\, &3\le k<d,
\\E_0, &k=d,
\end{cases}
$$
where $E_{d-k}$ is a parallelepiped in $\mathbb R ^k$ and
$$m_k (E_{d-k})= h^k m_{k-1}  (E_{d-k+1}).$$
Since $m_{k-1}  (E_{d-k+1})\le h^{\frac{(k-1)^2+(k-1)-2}{2}}
 \bigl(\phi ^{(d-1)}(t+h)-\phi ^{(d-1)}(t)\bigr)$
we also obtain
$$m_k (E_{d-k})\le
 h^{\frac{k^2+k-2}{2}}
 \bigl(\phi ^{(d-1)}(t+h)-\phi ^{(d-1)}(t)\bigr)
$$
and the claim is proved.

Finally, if we  apply the claim for $k=d$ and note that
 $\lambda (E_0)\ge h$, \eqref{hypothesis2} yields the conclusion of the
lemma.
\end{proof}

\begin{proof}[{\bf Proof of Lemma \ref{Lemma2}}]
We begin by noting an inequality for the Vandermonde determinant
\eqref{vanderm}, namely,  if $\delta >0$
and $t_1 <\dots <t_n$ then (with $u=(u_1,\dots, u_{n-1})$)
\begin{multline}\label{ineq9}
\int\nolimits_{t_1}^{t_2}\int\nolimits_{t_2}^{t_3}\cdots
\int\nolimits_{t_{n-1}}^{t_n}
V_{n-1}(u)(u_{n-1}-u_1 )^{\delta}du_{n-1}\cdots du_1
\\ \geq
C(n)\, V_n (t_1 ,\dots ,t_n )\, (t_n -t_1 )^{\delta}.
\end{multline}
To see \eqref{ineq9}, we observe that the left hand side is bounded below by
\begin{align*}
&\int\nolimits_{t_1}^{(t_1 +t_2)/2}\int\nolimits_{t_2}^{t_3}\cdots
\int\nolimits_{t_{n-2}}^{t_{n-1}}
\int\nolimits_{(t_{n-1}+t_n )/2}^{t_n}
V_{n-1}(u)(u_{n-1}-u_1 )^{\delta}du
\\
&\geq\Bigl(\frac{t_n -t_1}{2}\Bigr)^{\delta}
\int\nolimits_{t_1}^{(t_1 +t_2)/2}\int\nolimits_{t_2}^{t_3}\cdots
\int\nolimits_{t_{n-2}}^{t_{n-1}}
\int\nolimits_{(t_{n-1}+t_n )/2}^{t_n}
V_{n-1}(u)\, du.
\end{align*}
Now we also use  \eqref{vandintegration},
%
%
and together with the estimate
\begin{multline*}
\int\nolimits_{t_1}^{(t_1 +t_2 )/2}\int\nolimits_{(t_{n-1} +t_n )/2}^{t_n}
(u_{n-1}-u_1 )\prod\limits_{j=2}^{n-2}[(u_j -u_1 )(u_{n-1}-u_j)]\
du_{n-1}\, du_1
\\
 \geq
\frac{1}{4}\int\nolimits_{t_1}^{t_2}\int\nolimits_{t_{n-1}}^{t_n}
(u_{n-1}-u_1 )\prod\limits_{j=2}^{n-2}[(u_j -u_1 )(u_{n-1}-u_j)]\
du_{n-1}\, du_1,
\end{multline*}
this yields \eqref{ineq9}.

Now assume that the inequality \eqref{phiest}
holds if $a<s<t<b$ and let
$\cJ_k(t_1,\dots,t_k ;\phi^{(d-k)})$ be defined as in Lemma
\eqref{Jdsphilemma}, {\it  i.e.}, as the determinant of the
 $k\times k$ matrix with columns $(1,t_j ,\dots,
\tfrac{t_j^{k-2}}{(k-2)!},\phi^{(d-k+1}(t_j ))^T$.
We will show that if $2\le k\le d$ and
$a<t_1 <\cdots <t_k <b$, then
\begin{equation}\label{cJest}
\cJ _k (t_1 \dots ,t_k ;\phi ^{(d-k)})\ge c(k)\, c\, V_k (t_1 ,\dots ,t_k
)\, (t_k -t_1 )^ {\rho -1}.
\end{equation}

\noindent By choosing $\{ t_j \}$ to be a nondecreasing rearrangement of
$\{ s+h_j \}$, the case $k=d$ of \eqref{cJest}
will imply Lemma \ref{Lemma2}.
If $k=2$ then \eqref{cJest} follows immediately from \eqref{phiest}.
So, proceeding by induction, assume that \eqref{cJest} holds for $k-1$.
By \eqref{firstintegration}
\begin{multline*}
\cJ_k ((t_1 ,\dots ,t_k ;\phi^{(d-k)})\\=
\int_{t_1}^{t_2}\dots\int_{t_{k-1}}^{t_k }
\cJ_{k-1}(\sigma_1,\dots,\sigma_{k-1};\phi^{(d-k+1)})\,
d\sigma_{k-1}\cdots
d\sigma_{1}.
\end{multline*}
By our inductive assumption this exceeds
$$
c(k-1)\,c\,\int_{t_1}^{t_2}\dots\int_{t_{k-1}}^{t_k }
V_{k-1}(\sigma _1 ,\dots ,\sigma _{k-1} )(\sigma _{k-1}-\sigma _1 )^{\rho
-1}
d\sigma_{k-1}\cdots
d\sigma_{1}
$$ and so \eqref{ineq9} gives \eqref{cJest}, completing the proof of Lemma
\ref{Lemma2}.
\end{proof}

\section{Further results}\label{further}

In this section we gather some results  about Fourier restriction with
respect to Euclidean  arclength measure on curves,
mainly focusing on degenerate homogeneous curves.
For related  arguments see
\cite{sj},  \cite{so}, \cite{DM1}.

\subsection{Homogeneous curves}\label{homogen}
The following result follows by rescaling techniques from the result in
\cite{BOS} on nondegenerate curves (analogous to \eqref{moment}).
Let \begin{equation} \label{gahom}
\gamma(t)= (t^{a_1},t^{a_2}, \dots, t^{a_d})
\end{equation}
where $d\ge 3$, and $-\infty<a_1<a_2<\dots< a_d<\infty$, and
$a_i\neq 0$, $i=1,\dots, d$. We let $\cR$ be the Fourier restriction operator,
 setting
$\cR f(t)= \widehat f(\gamma(t))$. Let
$$D=a_1+a_2+\dots+a_d$$ be the ``homogeneous'' dimension and assume
$D> d(d+1)/2$.

\begin{proposition}
Let $p_d=\frac{d^2+d+2}{d^2+d}$ and $\gamma$ as in
\eqref{gahom}. Then
$\cR$ is of restricted weak type $(p_d, p_d'/D)$,
\begin{equation}\label{infty}
\big\|\cR f\big\|_{L^{p_d'/D,\infty}(dt)} \le C(a_1,\dots, a_n)
\|f\|_{L^{p_d,1}}.
\end{equation}
\end{proposition}

\begin{proof}
Define
\[ ({\cR}_k f) (t) = \widehat{f}(\gamma(t)) \chi_{I_k}(t) \]
where $I_k = [2^{-k-1}, 2^{-k}]$. We may use the nonisotropic dilations adapted
to the curve to rescale  the result in the
nondegenerate case (Theorem 1.1 in \cite{BOS}); we obtain
\begin{equation}\label{scaled}
\| {\cR}_k f\|_{L^{p_d}(dt)} \le C 2^{k[(D+1)(1-\frac 1{p_d})-1]} \|
f\|_{L^{p_d, 1}(\bbR^d)}.
\end{equation}
Let $D_0=d(d+1)/2$ and fix  $0<q_0 < p_d'/D$.
Since $1/q_0 > D/p_d' > D_0/p_d'= 1/p_d$, by H\"older's inequality
the last estimate implies
\begin{equation}\label{Holder}
\| {\cR}_k f\|_{L^{q_0}(dt)} \le C 2^{-k[\frac{1}{q_0} -D(1-\frac
1{p_d})]}
\| f\|_{L^{p_d,1}(\bbR^d)} .
\end{equation}
Since $(D+1)/(p_d')-1 > (D_0+1)/(p_d')-1 = 0$, an application of
Bourgain's interpolation lemma to \eqref{scaled} and
\eqref{Holder} gives the assertion.
\end{proof}

\subsection{An improvement}\label{improve}
For some very specific classes we can improve the second  Lorentz  exponent on the left hand side of \eqref{infty}.

We now suppose the stronger restricted strong type estimate
\begin{equation}\label{endpt}
\| {\cR} f\|_{L^{p_d}(w dt)} \le C \| f\|_{p_d, 1}
\end{equation}
where $w dt$ is affine arclength measure.
{\it Assume} that \begin{equation}\label{wassumpt}
1/w \in L^{s,\infty}(dt)\end{equation} for some $s\in (0, \infty)$. Define $q$ by
\begin{equation}\label{qass} \frac{1}{q} = \frac{1}{p_d} + \frac{1}{s p_d} .
\end{equation}
Then as in \cite{DM1} one can use
 the Lorentz space multiplication theorem (Theorem 4.5 in \cite{hu}),
and  it follows that
\begin{align*}
& \| \cR f \|_{L^{q, p_d}(dt)} = \| (\cR f) w^{1/p_d} \cdot w^{-1/p_d}
\|_{L^{q,
p_d}(dt)}
\\ &\le  C \| (\cR f) w^{1/p_d}\|_{L^{p_d}(dt)} \,
\| w^{-1/p_d}\|_{L^{s p_d,
\infty}(dt)}
=C\|w^{-1} \|_{L^{s,\infty}(dt)}^{1/p_d} \,\|{\cR} f\|_{L^{p_d}(w dt)} .
\end{align*}
Hence \eqref{endpt} and \eqref{wassumpt} imply that for $q$ as in \eqref{qass}
\begin{equation}\label{better}
\| {\cR} f\|_{L^{q, p_d}(dt)} \le C \| f\|_{p_d, 1}.
\end{equation}

\begin{corollary} Let  $\gamma(t)= (t, t^\alpha, t^{5\alpha-1})$ with
$\alpha >1$.  Then $\cR$ maps $L^{7/6,1}$ boundedly
to $L^{7/(6\alpha),7/6}$.
\end{corollary}
 \begin{proof} Note that  $D=6\alpha > 6= D_0$.
Also one computes
$w(t)=c(\alpha) t^{\alpha-1}$ with $c(\alpha)\neq 0$ so that
$w^{-1}\in L^{s,\infty}$ for $s=1/(\alpha-1)$.
By  Theorem 1.4
in \cite{BOS} it follows that
\eqref{endpt} holds  with $p_3=7/6$, so that the assertion follows.
\end{proof}


\subsection{$ {L}^{p}\to {L}^{{q}}$ bounds}\label{pq}
Finally, let us suppose that, instead of \eqref{endpt}, the
estimate
\begin{equation}\label{nonendpt} \| {\cR} f\|_{L^{Q}(w
dt)} \le C \| f\|_{p}
\end{equation}
holds for $1/p+ 1/(D_0 Q) =1$, and $1/w \in L^{s, \infty}(dt)$
with $1< p < p_d$ and some $s\in (0, \infty)$. Then an argument
similar to the one given above together with an interpolation
show that
\[
\| {\cR} f\|_{L^{q,p}(dt)} \le C\| f\|_{p} \] for $1<p<p_d$ and
\[ \frac{1}{p} + \frac{s}{(s+1)D_0 q} =1 .\]


\begin{thebibliography}{99}

\bibitem {baklee} J.-G. Bak, S. Lee, Estimates for an oscillatory
integral operator
related to restriction to space curves, \textit{Proc. Amer. Math. Soc.
}\textbf{132} (2004), 1393--1401.


\bibitem{BOS} J.-G. Bak, D. Oberlin,  A. Seeger,
Restriction of Fourier transforms to curves and related
oscillatory integrals, 
to appear in Amer. J. Math.

\bibitem{Bo1} J. Bourgain, Estimations de certaines
fonctions maximales,
\textit{C. R. Acad. Sci. Paris S\'er. I Math.}  \textbf{301}
(1985),
no. 10, 499--502.


\bibitem{CSWW} A. Carbery, A. Seeger, S. Wainger and J. Wright,
Classes of singular integral operators along variable lines,
\textit{J. Geom. Anal.} \textbf{9}  (1999),  no. 4, 583--605.


\bibitem{Ch} M. Christ, On the restriction of the Fourier transform to
curves:
endpoint results and the degenerate case,
\textit{Trans. Amer. Math. Soc.} \textbf{287} (1985), 223--238.


\bibitem{D} S.W. Drury, Restriction of Fourier transforms to curves,
\textit{Ann. Inst. Fourier}, \textbf{35} (1985), 117-123.

\bibitem{D2} \bysame, Degenerate  curves and harmonic analysis,
\textit{Math. Proc. Cambridge Philos.
Soc.} \textbf{108} (1990), 89-96.


\bibitem {DM1} S.W.  Drury, B. Marshall, Fourier restriction theorems for
curves
with affine and Euclidean arclengths,
\textit{Math. Proc. Cambridge Philos.
Soc.} \textbf{97} (1985), 111-125.

\bibitem {DM2}\bysame, Fourier restriction theorems for degenerate curves,
\textit{Math. Proc. Cambridge Philos. Soc.} \textbf{101} (1987), 541-553.


\bibitem{hu} R. Hunt, On $L(p,q)$ spaces, \textit{Enseignement Math.} \textbf{12} (1966), 249-276.



\bibitem{janson}
S. Janson,
On interpolation of multilinear operators. Function spaces and applications (Lund, 1986), 290--302,
Lecture Notes in Math., 1302, Springer, Berlin, 1988.



\bibitem{sj} P. Sj\"olin, Fourier multipliers and estimates of the Fourier
transform of measures carried by smooth curves in $R\sp{2}$,
\textit{Studia Math.}, \textbf{51} (1974), 169--182.


\bibitem{so} C.D. Sogge,  A sharp restriction theorem for degenerate curves in $R\sp 2$,
\textit{Amer. J. Math.} \textbf{109}  (1987),  no. 2, 223--228.


\end{thebibliography}
\end{document}